\pdfsuppressptexinfo=-1 \pdfinfoomitdate=1 \pdftrailerid{}
\documentclass[11pt]{amsart}
\usepackage[margin=1in]{geometry}
\usepackage{amssymb,amsmath,amsthm}
\usepackage[colorlinks=true,linkcolor=black,citecolor=black,urlcolor=black]{hyperref}
\AtBeginDocument{\hypersetup{pdftitle={A two-sided continued-fraction criterion for the Erdos-Moser equation and Erdos's 1949 problem},pdfauthor={Elham Mehdi-Nezhad and Amir M. Rahimi},pdfcreator={},pdfproducer={}}}

\newtheorem{theorem}{Theorem}
\newtheorem{proposition}[theorem]{Proposition}
\newtheorem{lemma}[theorem]{Lemma}
\newtheorem{corollary}[theorem]{Corollary}
\theoremstyle{remark}
\newtheorem{remark}[theorem]{Remark}
\newcommand{\gz}{\gamma_0}

\title[A two-sided continued-fraction criterion and Erd\H{o}s's 1949 problem]{A two-sided continued-fraction criterion for the\\ Erd\H{o}s--Moser equation and Erd\H{o}s's 1949 problem}
\author{Elham Mehdi-Nezhad}
\address{Department of Mathematics and Applied Mathematics, University of the Western Cape, Private Bag X17, Bellville 7535, Cape Town, South Africa}
\email{emehdinezhad@uwc.ac.za}
\author{Amir M. Rahimi}
\address{School of Mathematics, Institute for Research in Fundamental Sciences (IPM), P.O. Box 19395-5746, Tehran, Iran}
\email{amrahimi@ipm.ir}

\begin{document}
\begin{abstract}
We give a two-sided continued-fraction criterion for the Erd\H{o}s--Moser equation $1^k+2^k+\dots+(m-1)^k=m^k$. Writing $m=k/\log2+\frac32+\gamma_m/m$, Gallot, Moree and Zudilin showed that a solution yields a convergent of $\log2/(2N)$, for a known divisor $N$ of $k$, followed by a large partial quotient. We show that an explicit two-sided estimate for $\gamma_m$, proved in our earlier centre-of-mass work, restricts this partial quotient, for each value of $\gcd(k/N,2m-3)$, to two consecutive values, or to three in the exceptional case that an integer lies in an explicit interval of relative length $O(1/m)$. This excludes the first index left open by their computation and, without extending their continued-fraction computation, raises the lower bound for a nontrivial solution from $m>2.71\cdot10^{1\,667\,658\,416}$ to $m>6.26\cdot10^{1\,667\,658\,425}$.

The same criterion gives an effective form of the asymptotic results of Best and te Riele on Erd\H{o}s's 1949 problem on the largest term of $1^n+\dots+(m-1)^n$ (Advanced Problem 4347 in the Monthly). We determine all $50$ pairs $(m,n)$ with $m<10^{399}$ for which $(m-2)^n<\frac12(m-1)^n$ and $1^n+\dots+(m-1)^n>m^n$; in particular, the smallest example found by Best and te Riele in 1976 is the smallest such pair, a question they had left open. For $n<0.69\cdot10^{399}$ we also determine the $51$ elements of O'Bryant's exceptional set $\mathcal A_2$ for $\lfloor1/(2^{1/n}-1)\rfloor$ in this range and all $101$ terms of OEIS A129935 below this bound, one of them not previously listed, and settle several conjectures recorded in the OEIS. Finally, using only the upper bound in our estimate for $\gamma_m$ and the positivity of $\gamma_m$, we show that the second part of Erd\H{o}s's problem, the existence of infinitely many such pairs, is equivalent to the infinitude of $\mathcal A_2$ and to the existence of infinitely many even-index convergents $p_j/q_j$ of $\log2$ with $p_j$ even, $q_j$ odd and
\[
q_j(q_j\log2-p_j)<\Theta(q_j),
\]
where
\[
\Theta(q)=\frac{4\log2\;h_{(q+3)/2}}{1+3/q},\qquad h_m=m\Big(m-\frac32+\frac1{\log(1-\frac1{m-1})}\Big),
\]
so that $\Theta(q)=\frac{\log2}{3}+O(q^{-2})$. In particular, the condition no longer involves $\gamma_m$.
\end{abstract}

\subjclass[2020]{Primary 11D61; Secondary 11A55, 11J70, 11Y50}
\keywords{Erd\H{o}s--Moser equation, continued fractions, $\log 2$}
\maketitle

\section{Introduction}

Erd\H{o}s conjectured that
\begin{equation}\label{eq:EM}
1^k+2^k+\dots+(m-1)^k=m^k
\end{equation}
has no solution in positive integers other than $1+2=3$. Moser~\cite{Moser} proved that any further solution has $m>10^{10^6}$; see~\cite{Moree} for the history. Let $L=\log 2$. Gallot, Moree and Zudilin~\cite{GMZ} introduced continued fractions into the problem. They showed that $2k/(2m-3)$ is a convergent of $L$ of even index; combining this with divisibility properties of $k$ and $2m-3$ and with the continued fraction of $L/(2N)$ for a suitable integer $N$, they obtained $m>2.7139\cdot10^{1\,667\,658\,416}$, which to our knowledge is still the best bound known.

Let $\kappa_m$ be the real number with $1^{\kappa_m}+\dots+(m-1)^{\kappa_m}=m^{\kappa_m}$, and define $\gamma_m$ by
\begin{equation}\label{eq:gamma}
\kappa_m=L\Big(m-\frac32-\frac{\gamma_m}{m}\Big),\qquad \gz=\frac{25}{12}-3L=0.0038917\ldots
\end{equation}
The solutions of~\eqref{eq:EM} are the $m$ for which $\kappa_m=k$ is an integer. The method of~\cite{GMZ} uses the analytic information about $\gamma_m$ only through the bounds $0<\gamma_m<0.004$; the upper bound yields the inequality $a_{j+1}\ge180N-2$ for the partial quotient following the relevant convergent (condition~(b) of~\cite[Theorem~2]{GMZ}). The authors show that $\gamma_m=\gz+O(1/m)$ and remark that they have no general strategy for making further terms of the expansion explicit. We address this point: a lower bound for $\gamma_m$ gives an upper bound for $a_{j+1}$ as well, and an explicit two-sided bound restricts $a_{j+1}$, for each value of $g$, to two consecutive values, or to three in the exceptional case that an integer lies in an explicit interval of relative length $O(1/m)$ (Proposition~\ref{prop:pin} and Remark~\ref{rem:heur}). (A one-sided form of this identity was used in~\cite[Theorem~8.2(b)]{MR} for the analogue of~\eqref{eq:EM} in arithmetic progressions.) The bound we use is our estimate~\cite[Theorem~5.5]{MR}
\begin{equation}\label{eq:thm55}
-\frac{87}{m}<\gamma_m-\gz<\frac{116}{m}\qquad(m\ge10^5),
\end{equation}
obtained by comparing the power sum with the moments of a geometric mass distribution (the mass-system viewpoint for sums of $k$-th powers goes back to~\cite{R09,R17}). Before~\eqref{eq:thm55}, the only explicit lower bound for $\gamma_m$ we know of was $\gamma_m>0$. The lower bound in~\eqref{eq:thm55} is essential: the upper bound $\gamma_m<0.004$ only says that $a_{j+1}$ is large, which happens at many indices, whereas the two bounds together leave two possible values of $a_{j+1}$ (three in the exceptional case described in Remark~\ref{rem:heur}) for each value of $g=\gcd(k/N,2m-3)$. Heuristically, for a number whose partial quotients follow the Gauss--Kuzmin law, the probability that a given even index passes the test drops from about $1/(125N)$ to about $1/(1.2\cdot10^4N^2)$.

The only analytic estimates from~\cite{MR} used here are~\eqref{eq:thm55} and the bounds $0<\gamma_m<1/(8L)$ of~\cite[Theorem~6.6]{MR}. The main new ingredient is the two-sided criterion (Proposition~\ref{prop:pin} and Corollary~\ref{cor:crit}), which replaces the one-sided use of the identity in~\cite[Theorem~8.2(b)]{MR}; Theorems~\ref{thm:main} and~\ref{thm:M}, Proposition~\ref{prop:inf}, Corollary~\ref{cor:inf} and the deductions in Remarks~\ref{rem:oeis} and~\ref{rem:A2} are new.

\begin{theorem}\label{thm:main}
Every solution of~\eqref{eq:EM} with $k\ge2$ satisfies $m>6.2612\cdot10^{1\,667\,658\,425}$.
\end{theorem}

This improves the bound of~\cite{GMZ} by a factor of about $2.3\cdot10^9$ without new computation: we exclude the first index that their criteria leave open. With the new criterion the bound is limited only by the length of the continued-fraction computation (Section~\ref{sec:comp}).

With $\alpha=L$ and no auxiliary integer $N$, the identity applies to an older problem of Erd\H{o}s~\cite{Er49} on the largest term of $1^n+\dots+(m-1)^n$ (Section~\ref{sec:erdos}). Best and te Riele~\cite{BtR} settled its first part and found $33$ pairs $(m,n)$ for the second, but could not show that no such pair is smaller than their smallest example; they observed that this would need effective versions of their asymptotic results. We give such a version. Theorem~\ref{thm:M} determines the $50$ pairs with $m<10^{399}$ and shows that the example of Best and te Riele is the smallest. This minimality also follows, through the correspondence in Remark~\ref{rem:A2}, from O'Bryant's determination~\cite[\S6]{OB} of the elements of the exceptional set $\mathcal A_2$ of $\lfloor1/(2^{1/n}-1)\rfloor$ up to $777\,451\,915\,729\,368$; Remark~\ref{rem:A2}, using the proof of Theorem~\ref{thm:M}, extends this determination to $n<0.69\cdot10^{399}$ (these values belong to the larger set OEIS A129935, whose terms below this bound we also determine: the $100$ listed there and one more). Remarks~\ref{rem:oeis} and~\ref{rem:A2} also show that several conjectures in the OEIS are decided in this range, among them one of Greathouse, which fails at exactly $51$ arguments $N<10^{399}-2$. Finally, Proposition~\ref{prop:inf} shows that the upper bound in~\eqref{eq:thm55} removes $\gamma_m$ from the question of infinitude: part~(ii) of Erd\H{o}s's problem holds if and only if $\mathcal A_2$ is infinite, and if and only if infinitely many convergents $p_j/q_j$ of $L$ of even index with $p_j$ even and $q_j$ odd satisfy $q_j(q_jL-p_j)<\Theta(q_j)$, where $\Theta(q_j)=\frac L3+O(q_j^{-2})$ is explicit (Corollary~\ref{cor:inf}). The last equivalence sharpens criteria of O'Bryant~\cite[Lemmas~11 and~12]{OB}; the equivalence with part~(ii) is new.

\section{The partial quotient at a solution}

\begin{lemma}[{\cite[Ch.~X]{HW}}]\label{lem:cf}
Let $\alpha=[a_0;a_1,\dots]$ be irrational, with convergents $p_i/q_i$ and complete quotients $\alpha_i=[a_i;a_{i+1},\dots]$. Then
$q_j^2\,|\alpha-p_j/q_j|=\big(\alpha_{j+1}+q_{j-1}/q_j\big)^{-1}$.
\end{lemma}

\begin{proposition}\label{prop:pin}
Let $(m,k)$ be a solution of~\eqref{eq:EM} with $k\ge2$, let $N\mid k$, $k=Nk_1$, $\alpha=L/(2N)$ and $g=\gcd(k_1,2m-3)$, which is odd since $2m-3$ is. Then $k_1/(2m-3)=p_j/q_j$ in lowest terms is a convergent of $\alpha$ of even index $j$, $2m-3=gq_j$, and
\begin{equation}\label{eq:pin}
\alpha_{j+1}+\frac{q_{j-1}}{q_j}=B_mg^2\Big(1+\frac{3}{gq_j}\Big),\qquad B_m=\frac{N}{2\gamma_mL}.
\end{equation}
In particular, for $j\ge2$,
\begin{equation}\label{eq:window}
B_mg^2-2<a_{j+1}<B_mg^2\Big(1+\frac{3}{gq_j}\Big).
\end{equation}
\end{proposition}

\begin{proof}
By~\eqref{eq:gamma} with $\kappa_m=k$ we have
$L-\frac{2k}{2m-3}=\frac{2\gamma_mL}{m(2m-3)}$. Dividing by $2N$ and writing $q=2m-3$, $m=(q+3)/2$, $q=gq_j$, $k_1=gp_j$, we get
\[
\alpha-\frac{p_j}{q_j}=\frac{2\gamma_mL}{Ng^2q_j(q_j+3/g)} .
\]
By Moser's bound and $0<\gamma_m<0.004$~\cite{GMZ}, the right side is positive and less than $1/(2q_j^2)$. By Legendre's theorem~\cite[Theorem~184]{HW}, $p_j/q_j$ is a convergent; since the convergents of even index are exactly those below $\alpha$, its index is even. Lemma~\ref{lem:cf} now gives~\eqref{eq:pin}. Finally,~\eqref{eq:window} follows from $a_{j+1}\le\alpha_{j+1}<a_{j+1}+1$ and $0<q_{j-1}/q_j<1$ for $j\ge2$.
\end{proof}

\begin{remark}
Identity~\eqref{eq:pin} in fact determines $\gamma_m$:
\[
\gamma_m=\frac{N\,(1+3/(gq_j))}{2Lg^2\,\big(\alpha_{j+1}+q_{j-1}/q_j\big)}.
\]
So the test can be sharpened further: this value, computed from $\alpha_{j+1}$ and $q_{j-1}/q_j$, must lie in $(\gz-87/m,\gz+116/m)$.
\end{remark}

\begin{corollary}[criterion]\label{cor:crit}
Suppose that $N$ divides $k$ for every solution of~\eqref{eq:EM} with $k\ge2$. Let $j\ge2$ be even with $q_j\ge2\cdot10^5$. Suppose that for every odd $g\ge1$ and every $\gamma\in(\gz-87/m,\gz+116/m)$, where $m=(gq_j+3)/2$, the integer $a_{j+1}$ does not satisfy~\eqref{eq:window} with $B_m$ replaced by $N/(2\gamma L)$. Then no solution of~\eqref{eq:EM} corresponds to the convergent $p_j/q_j$.
\end{corollary}

\begin{proof}
Combine Proposition~\ref{prop:pin} with~\eqref{eq:thm55}.
\end{proof}

In practice, for a given even index $j$ one checks the finitely many odd $g$ with $Bg^2-2\le a_{j+1}+1$, and for each of them whether $a_{j+1}$ lies in the interval~\eqref{eq:window}, using the extreme values $\gz-87/m$ and $\gz+116/m$ of $\gamma$. The upper bound for $\gamma_m$ enters only through the left end of~\eqref{eq:window}, as in~\cite{GMZ}; the lower bound gives the right end, and it is the right end that makes the set of admissible $a_{j+1}$ finite for each $g$.

\begin{remark}\label{rem:heur}
Put $B=N/(2\gz L)=185.351\ldots N$. Let $X=B_mg^2(1+3/(gq_j))$ be the right side of~\eqref{eq:pin}. Since $a_{j+1}\le\alpha_{j+1}<a_{j+1}+1$ and $0<q_{j-1}/q_j<1$, we have $a_{j+1}<X<a_{j+1}+2$, so for the exact value of $\gamma_m$ there are at most two possible values of $a_{j+1}$, namely $\lceil X\rceil-2$ and $\lceil X\rceil-1$. By~\eqref{eq:thm55}, $X$ is known to within a relative error below $203/(\gz m)+O(m^{-2})<6\cdot10^4/m$, so $a_{j+1}$ is restricted to two consecutive integers unless an integer lies in this short interval around $X$, in which case to three. For $m$ of the size relevant here the relative error is below $10^{-1.6\cdot10^9}$ and $X$ is close to $Bg^2$. Condition~(b) of~\cite[Theorem~2]{GMZ}, $a_{j+1}\ge180N-2$, is the one-sided version of this. Heuristically, over the roughly $1.6\cdot10^9$ even indices up to the index $j^*$ below, about $1.7$ indices are expected to pass condition~(b) (Table~1 of~\cite{GMZ} records two indices up to $j^*$ passing (a), (b) and (c)), against about $2.3\cdot10^{-9}$ for the window~\eqref{eq:window}.
\end{remark}

\section{Proof of Theorem~\ref{thm:main}}

Let $N=2^8\cdot3^5\cdot5^3=7\,776\,000$. By Moree, te Riele and Urbanowicz~\cite[p.~814]{MtRU}, combined with a result of Kellner~\cite{Kellner} (see~\cite[\S1]{GMZ}), the number $N_2=2^8\cdot3^5\cdot5^4\cdot7^3\cdots$ divides $k$ for every solution with $k\ge2$; in particular $N\mid k$. By~\cite[Theorem~2 and Table~1]{GMZ}, the index $j^*=3\,236\,170\,820$ is the smallest even index of the continued fraction of $\alpha=L/(2N)$ satisfying conditions (b)--(d) of~\cite[Theorem~2]{GMZ}. By the proof of~\cite[Theorem~2]{GMZ}, the convergent attached to any solution satisfies (a)--(d), whatever the value of $g=\gcd(k_1,2m-3)$: (a) holds by Proposition~\ref{prop:pin}; (b) holds since~\eqref{eq:window} gives $a_{j+1}>B_mg^2-2\ge B_m-2>180N-2$, because $\gamma_m<0.004$; (c) holds since $q_j\mid 2m-3$ and $(2m-3,6)=1$ (that is, $3\nmid m$; see the proof of~\cite[Theorem~2]{GMZ}). For (d), write $\nu_p$ for the $p$-adic valuation. If $p\in\mathcal P(N)$ divides $q_j$, then $p\mid 2m-3$ and $p-1\mid k$, and~\cite[Lemma~2]{GMZ} gives $\nu_p(2m-3)=\nu_p(3^{p-1}-1)+\nu_p(k)+1$; here $\mathcal P(N)$ is the set of primes $p$ with $p-1\mid N$ or with $3$ a primitive root modulo $p$, as in~\cite{GMZ}. Since $q_j=(2m-3)/g$ with $g=\gcd(k_1,2m-3)$ and $\nu_p(q_j)\ge1$, we have $\nu_p(g)=\nu_p(k_1)$, so $\nu_p(q_j)=\nu_p(2m-3)-\nu_p(k_1)$; with $\nu_p(k)=\nu_p(N)+\nu_p(k_1)$ this gives $\nu_p(q_j)=\nu_p(3^{p-1}-1)+\nu_p(N)+1$, which is (d). Hence every solution corresponds to an even index $j\ge j^*$. Table~1 of~\cite{GMZ} gives
\[
a_{j^*+1}=2\,307\,115\,390,\qquad q_{j^*}>5.427815\cdot10^{1\,667\,658\,416}.
\]

It remains to exclude $j^*$ itself; this is the only step that uses the lower bound in~\eqref{eq:thm55}. Suppose that a solution corresponds to $j^*$. Since $m>10^5$, by~\eqref{eq:thm55} we have $\gz-0.00087<\gamma_m<\gz+0.00116$. Since $g$ is odd, \eqref{eq:window} gives
\[
a_{j^*+1}<B_m\Big(1+\frac3{q_{j^*}}\Big)<1.857\cdot10^9\quad(g=1),\qquad a_{j^*+1}>9B_m-2>9.9\cdot10^9\quad(g\ge3).
\]
Both contradict $a_{j^*+1}=2\,307\,115\,390$. (For comparison, $B=N/(2\gz L)=1\,441\,289\,467.265\ldots$ and $a_{j^*+1}/B=1.6007\ldots$.) Hence the index $j_s$ attached to a solution, which is even and at least $j^*$, satisfies $j_s\ge j^*+2$, and
\[
m=\frac{gq_{j_s}+3}{2}>\frac{q_{j^*+2}}{2}\ge\frac{q_{j^*+1}+q_{j^*}}{2}>\frac{(a_{j^*+1}+1)\,q_{j^*}}{2}>6.2612\cdot10^{1\,667\,658\,425}.\qquad\qed
\]

\begin{remark}
For $g=1$ the lower bound $\gamma_m>0.00244$ suffices, and an explicit form of the error term in~\cite[Theorem~1]{GMZ} would give this equally well; for $g\ge3$ the bound $\gamma_m<0.004$ of~\cite{GMZ} suffices. The same argument excludes the index $j=2\,015\,279\,170$ of Table~1 of~\cite{GMZ} ($a_{j+1}/B=3.045\ldots$) independently of their condition~(d). As a consistency check, the row $N=2^8\cdot3^5\cdot5^4$ of the same table ($j=3\,236\,257\,942$, $a_{j+1}=11\,535\,576\,954$, the same $q_j$) gives $a_{j+1}/(5B)=1.6007\ldots$, as expected.
\end{remark}

\section{Further computation}\label{sec:comp}

This section is not used in the proof of Theorem~\ref{thm:main}; it estimates the potential reach of the method with longer continued fractions, and no bound stated here has been computed or certified. By the heuristic in Remark~\ref{rem:heur}, the criterion is not expected to leave spurious candidates. The method should therefore certify $m>q_r/2$, where $r$ is the last index of the continued fraction of $L/(2N)$ that can be computed reliably. By the theorems of Lochs~\cite{Lochs} and L\'evy (see~\cite[\S4]{GMZ} and~\cite{Kh}), $d$ decimal digits of $L$ give $r\approx0.970\,d$ and $\log_{10}q_r\approx0.515\,r$. Since the product of these two constants is $1/2$, the result is a bound of roughly $10^{d/2}$. If the continued fraction computed in~\cite{GMZ} beyond $j^*$ is still available, Theorem~\ref{thm:main} extends to the last index computed there, unless some further index passes the test of Corollary~\ref{cor:crit}. Gallot, Moree and Zudilin note that $10^{10^{10}}$ is within reach with enough memory~\cite{GMZ}. Since $L$ is now known to $3.3\cdot10^{12}$ decimal digits~\cite{YC}, a bound of order $10^{1.65\cdot10^{12}}$ is in principle accessible.

\section{Erd\H{o}s's problem on the largest term}\label{sec:erdos}

Let $S_n(m)=1^n+\dots+(m-1)^n$. In Advanced Problem 4347 of the Monthly~\cite{Er49}, Erd\H{o}s conjectured, in the formulation of~\cite{BtR}:
\begin{enumerate}
\item[(i)] $S_n(m)\ge m^n$ whenever $(m-2)^n\ge\frac12(m-1)^n$;
\item[(ii)] $S_n(m)>m^n$ for infinitely many pairs with $(m-2)^n<\frac12(m-1)^n$.
\end{enumerate}
Van de Lune and te Riele~\cite{vdLtR} (see also~\cite{vdL}) proved (i) for a set of $n$ of density $1$, and Best and te Riele~\cite{BtR} proved (i) in full. Let $\mathcal M$ be the set of pairs in (ii). Best and te Riele showed that $\#\{(m,n)\in\mathcal M:m\le x\}=O(\log x)$ and found $33$ elements from the continued fraction of $2/\log2$, the smallest having $m=1\,121\,626\,023\,352\,385$ and $n=777\,451\,915\,729\,368$.

Put $\ell_m=L/(-\log(1-\frac1{m-1}))$, so that $(m-2)^n<\frac12(m-1)^n$ if and only if $n>\ell_m$. Since $n\mapsto S_n(m)/m^n$ is decreasing, $S_n(m)>m^n$ if and only if $n<\kappa_m$. Hence
\[
(m,n)\in\mathcal M\iff \ell_m<n<\kappa_m .
\]
Write $n=L(m-\frac32-\frac\delta m)$, which defines $\delta=\delta(m,n)$; thus $\delta$ plays for the integer $n$ the role that $\gamma_m$ plays for $\kappa_m$ in~\eqref{eq:gamma}, and comparing $n$ with $\kappa_m$ amounts to comparing $\delta$ with $\gamma_m$. Then $n<\kappa_m$ if and only if $\delta>\gamma_m$, and $n>\ell_m$ if and only if $\delta<h_m$, where
\[
h_m=m\Big(m-\frac32+\frac1{\log(1-\frac1{m-1})}\Big)
\]
decreases from $h_3=0.1719\ldots$ to $\frac1{12}$. Indeed, with the Gregory coefficients $G_i=\int_0^1\binom{u}{i}\,du$, which satisfy $(-1)^{i-1}G_i>0$ (for $0<u<1$ the product $u(u-1)\cdots(u-i+1)$ has exactly $i-1$ negative factors), one has $-1/\log(1-x)=1/x-\sum_{i\ge1}|G_i|x^{i-1}$ with $|G_1|=\frac12$, $|G_2|=\frac1{12}$, hence
\[
h_m=\frac{m}{m-1}\sum_{i\ge2}|G_i|\,(m-1)^{-(i-2)},
\]
a sum of positive terms decreasing in $m$, with limit $|G_2|=\frac1{12}$. So $\mathcal M$ consists of the pairs with $\gamma_m<\delta<h_m$, a window of length about $\frac1{12}-\gz=3L-2$ in $\delta$.

\begin{proposition}\label{prop:M}
Let $(m,n)\in\mathcal M$, $q=2m-3$ and $t=\gcd(2n,q)$, which is odd. Then $2n/q=p_j/q_j$ in lowest terms is a convergent of $L$ of even index, with $p_j$ even, $q_j$ odd, $q=tq_j$, and, with $\alpha_{j+1}$ the complete quotient of $L$,
\begin{equation}\label{eq:M}
\alpha_{j+1}+\frac{q_{j-1}}{q_j}=\frac{t^2}{4L\delta}\Big(1+\frac3{tq_j}\Big),\qquad \gamma_m<\delta<h_m .
\end{equation}
Conversely, every even $j$ with $p_j$ even and $q_j$ odd, and every positive odd $t$ for which the value of $\delta$ given by~\eqref{eq:M}, with $m=(tq_j+3)/2$, lies in $(\gamma_m,h_m)$, gives the element $(m,tp_j/2)$ of $\mathcal M$.
\end{proposition}

\begin{proof}
From the definition of $\delta$, $qL-2n=2L\delta/m=4L\delta/(q+3)$, so $0<q(qL-2n)<4Lh_3<\frac12$; positivity uses $\delta>\gamma_m>0$, valid for all $m\ge3$~\cite[Theorem~6.6]{MR}. By Legendre's theorem $2n/q$ is a convergent of even index, and Lemma~\ref{lem:cf} with $\alpha=L$ gives~\eqref{eq:M} as in the proof of Proposition~\ref{prop:pin}. The converse is the same computation read backwards.
\end{proof}

Formally,~\eqref{eq:M} is~\eqref{eq:pin} with $N=\frac12$, the gcd $g$ renamed $t$, and $\gamma_m$ replaced by $\delta$. Since $\delta\in(\gamma_m,h_m)$, with $\gamma_m\approx\gz$ and $h_m\to\frac1{12}$, for $t=1$ the left side of~\eqref{eq:M} lies essentially in $(3/L,\,1/(4L\gz))=(4.33,\,92.7)$, so $3\le a_{j+1}\le92$; for $m\ge10^5$, $a_{j+1}>4.33\,t^2-2$ is needed. The value $a_{j+1}=3$ does occur ($j=180$ and $j=422$). The upper end $92.7$ restricts only the value $t=1$: the windows $(4.33\,t^2,92.7\,t^2)$ for consecutive odd $t$ overlap, and Proposition~\ref{prop:inf} below shows that, once all odd $t$ are allowed, only the lower end matters.

\begin{theorem}\label{thm:M}
There are exactly $50$ pairs $(m,n)\in\mathcal M$ with $m<10^{399}$. The smallest has $m=1\,121\,626\,023\,352\,385$ and $n=777\,451\,915\,729\,368$; in particular the smallest example of Best and te Riele is the smallest element of $\mathcal M$.
\end{theorem}

\begin{proof}
Let $\tilde L$ be $L$ truncated to $3000$ decimals, so that $|L-\tilde L|<10^{-3000}$. The continued fractions of the rationals $\tilde L\pm10^{-2999}$ agree in their first $2925$ partial quotients, hence so does every real number between them, in particular $L$; this covers the $809$ convergents with $q_j<10^{400}$ and the partial quotients $a_{j+1}$ following them. By Proposition~\ref{prop:M} every element with $m<10^{399}$ comes from one of these convergents (since $q_j\le2m-3<10^{400}$) and an odd $t$ with $t^2<4Lh_3(a_{j+1}+2)$, where $4Lh_3=0.4766\ldots$; indeed, by~\eqref{eq:M} and $\delta<h_m\le h_3$, we have $t^2/(4Lh_3)<\alpha_{j+1}+q_{j-1}/q_j<a_{j+1}+2$. We tested the $149$ candidates $(m,n)=((tq_j+3)/2,\,tp_j/2)$ with $m<10^{399}$ and odd $t$ satisfying the weaker condition $t^2<\frac12(a_{j+1}+2)$, as required in Remark~\ref{rem:oeis}; for each we computed $\delta=m(m-\frac32-n/L)$ and $h_m$ in decimal arithmetic with $3200$ significant digits (the Python \texttt{decimal} module, in which every operation, including $\log$ and $\exp$, is correctly rounded, so each step has relative error at most $u=10^{-3199}$; $m$ and $n$ are represented exactly). Accuracy is lost only by cancellation in $m-\frac32-n/L$ and in $\log(1-\frac1{m-1})$, and the loss is bounded as follows. The computed $n/L<10^{399}$ has absolute error below $3u\cdot10^{399}<10^{-2799}$, the subtraction from $m-\frac32$ is exact, and multiplication by $m<10^{399}$ gives an error below $10^{-2399}$ in $\delta$. The computed $1-\frac1{m-1}\ge\frac12$ has absolute error at most $2u$, so its logarithm, of modulus at least $1/m$, has relative error at most $6um$; the reciprocal, of modulus at most $m$, has absolute error at most $7um^2$, and $h_m$ has absolute error at most $8um^3<10^{-2000}$. These bounds are smaller by more than $1900$ orders of magnitude than the margins stated next (for the $51$ values below one even has $|\delta-\gz|-116/m>2.9\cdot10^{-3}$), so every comparison below is certified. Three candidates have $m<10^5$ ($j=2,6,10$, with $m=3,140,1147$); all have $\delta\ge h_m$, i.e.\ $n\le\ell_m$, so they are not in $\mathcal M$ ($m=3$, $n=1$ is the equality case $(m-2)^n=\frac12(m-1)^n$). All candidates with $m\ge10^5$ satisfy $|\delta-h_m|>2\cdot10^{-4}$, and exactly $51$ of them satisfy $\delta<h_m$. Each of these $51$ satisfies $|\delta-\gz|\,m>2.2\cdot10^{13}$, so~\eqref{eq:thm55} decides whether $\delta>\gamma_m$. It holds for $50$ of them and fails for $(j,t)=(658,1)$, where $\delta=0.00094<\gz-87/m$.
\end{proof}

The first elements, by the index $j$ of the convergent of $L$, are $(j,t,a_{j+1})=(40,1,14)$, $(44,1,4)$, $(46,1,18)$, $(60,1,4)$, $(84,1,4)$, with $n$ of $15$, $18$, $19$, $26$ and $36$ digits; the index $j=658$, with $a_{j+1}=381$, gives four elements ($t=3,5,7,9$), while $t=1$ fails there because $\delta<\gamma_m$. Table~2 of~\cite{BtR} consists of the $33$ elements with $t=1$ and $j\le606$; the two elements with $(j,t)=(108,3)$ and $(494,3)$ lie in the same range but are missing there, since only fractions $2n/(2m-3)$ in lowest terms ($t=1$) were considered. (In~\cite{BtR} the convergents of $2/\log2$ are used, so the indices there differ from ours.)

\begin{remark}[Conjectures in the OEIS]\label{rem:oeis}
\emph{A072633.} The candidate $(j,t)=(658,1)$ excluded in the proof has $0<\delta<\gamma_m$; this is the exceptional exponent with $330$ digits noted after~\cite[Theorem~8.2]{MR}. Let $n\ge1$. Since $S_n(m)\ge m^n$ if and only if $n\le\kappa_m$, and $0<\gamma_M<1/(8L)$ for all $M\ge3$~\cite[Theorem~6.6]{MR}, the least $m'$ with $1^n+\dots+m'^n\ge(m'+1)^n$ is $\lceil n/L+\frac12\rceil$, except when $0<\delta(M,n)<\gamma_M$ for $M=\lceil n/L+\frac32\rceil$, in which case it is one larger (if $\delta=\gamma_M$, that is $n=\kappa_M$, equality holds and the value is $\lceil n/L+\frac12\rceil$). The pair $(j,t)=(658,1)$ is of this kind, with $n=p_{658}/2$, so the conjecture $a(n)=1+\operatorname{round}(n/\log2)$ recorded by H.~Kociemba for sequence A072633 of the OEIS~\cite{OEIS} (checked there for $n\le700$) fails at this $n$. This counterexample was first announced on 22 July 2026, with a verification by interval arithmetic, in~\cite{asanso}; the argument here does not depend on it. There is no other with $n<0.69\cdot10^{399}$: a counterexample needs $0<\delta<\gamma_M<1/(8L)$, which forces $t^2<\frac12(a_{j+1}+2)$; the candidates tested in the proof of Theorem~\ref{thm:M} include all of these, and a candidate with $\delta\ge h_m$ cannot have $\delta<\gamma_m$, since $\gamma_m<0.0051<\frac1{12}<h_m$ for $m\ge10^5$. Hence, by that computation together with a direct check of $S_n(m)\ge m^n$ for $(m,n)=(3,1)$ (equality), $(140,96)$ and $(1147,794)$, it is the only counterexample with $n<0.69\cdot10^{399}$.

\emph{A332101, A078607 and A332097.} Let $n\ge2$. Let $b(n)$ be the least $m$ with $m^n\le S_n(m)$ (sequence A332101; $b(n)$ is A072633$(n)+1$), and let $c(n)$ be the least $x$ with $2x^n>(x+1)^n$ (sequence A078607), so that $c(n)=\lceil1/(2^{1/n}-1)\rceil$ and $c(n)+2$ is the least $m$ with $n<\ell_m$; here $n\ne\ell_m$ for every $m$, since $2^{1/n}$ is irrational. By part~(i) of the conjecture, proved in~\cite{BtR}, $b(n)\le c(n)+2$; and since $\kappa_m<L(m-\frac32)<\ell_{m+1}$, also $b(n)\ge c(n)+1$. Thus $b(n)=c(n)+1$ exactly when $\ell_m<n\le\kappa_m$ for $m=b(n)$, that is, when $(m,n)\in\mathcal M$ or $(m,n)$ is a solution of~\eqref{eq:EM}. By Moser's bound there are no such solutions with $n<0.69\cdot10^{399}$, so in this range $b(n)=c(n)+2$ fails exactly at the $50$ values of $n$ from Theorem~\ref{thm:M}, the first being $777\,451\,915\,729\,368$. The formula $b(n)=c(n)+2$, stated in A332101 to hold for almost all $n$ with that $n$ as a possible exception, and the equivalent comment in A078607 that $c(n)+2$ is the least $s$ with $s^n<S_n(s)$, hold for almost all $n$ by~(i) and the bound $O(\log x)$ of~\cite{BtR}; the above shows that $777\,451\,915\,729\,368$ is indeed an exception and determines all exceptions in this range. (That the two formulas of A332101 disagree at this $n$ was noted there by T.~Sungkawichai in 2026.) For $E(s)=s^n-S_n(s)$ one has $E(s+1)-E(s)=(s+1)^n-2s^n$, so the maximum of $E$ is attained exactly at $s=c(n)$; hence the remark in A332097 that this maximum is almost always attained at $b(n)-2$ holds exactly when $b(n)=c(n)+2$. The conjecture $b(n)=\operatorname{round}(n/L)+2$, also recorded in A332101, is that of A072633 and fails at the same $330$-digit $n$. Finally, $c(n)+2$ differs from $\operatorname{round}(n/L)+2=\lceil n/L+\frac32\rceil$ exactly when $0<\delta(M,n)<h_M$ for $M=\lceil n/L+\frac32\rceil$; the argument of Proposition~\ref{prop:M} applies to such pairs, and by the proof of Theorem~\ref{thm:M} there are exactly $51$ of them with $n<0.69\cdot10^{399}$, namely the $50$ above and the $330$-digit $n$ with $(j,t)=(658,1)$. So $c(n)$ is the nearest integer to $n/\log2$ for all other $n$ in this range, the exception noted in A078607 being the first.
\end{remark}

\begin{remark}[The exceptional set $\mathcal A_2$, A129935 and A175406]\label{rem:A2}
This remark relates the $51$ values of $n$ with $0<\delta<h_m$ found in the proof of Theorem~\ref{thm:M} to O'Bryant's set $\mathcal A_2$, to the sequence A129935, which strictly contains it, and to A175406. Only O'Bryant's bound $777\,451\,915\,729\,368$ and the definitions of the sequences are taken from earlier work; the determinations beyond O'Bryant's bound are new.

\emph{The set $\mathcal A_2$.} These $51$ values of $n$ form the set $\mathcal A_2\cap[2,0.69\cdot10^{399})$, where $\mathcal A_\theta$ is the set of $n$ with $\lfloor1/(\theta^{1/n}-1)\rfloor\ne\lfloor n/\log\theta-\frac12\rfloor$ studied by O'Bryant~\cite{OB}; for $n\ge2$ the condition for $\theta=2$ reads $c(n)\ne\operatorname{round}(n/L)$. O'Bryant~\cite[\S6]{OB} showed that $\mathcal A_2$ has no element strictly between $1$ and $777\,451\,915\,729\,368$, and that this number lies in $\mathcal A_2$ (the discrepancy at this $n$ is recorded in A078608, where it is credited to D.~Hickerson, 2002). Since every $(m,n)\in\mathcal M$ gives an element $n$ of $\mathcal A_2$, his computation together with this correspondence also proves the minimality statement of Theorem~\ref{thm:M}. All of the $51$ values except the largest, the $393$-digit value from $j=796$, appear in sequence A129935 (D.~Hickerson, M.~Alekseyev, R.~Gerbicz; the mechanism was explained by D.~Applegate), which is defined by $\lceil2/(2^{1/n}-1)\rceil\ne\lfloor2n/L\rfloor$; this condition also holds when $n/L$ lies just below an integer, so A129935 strictly contains $\mathcal A_2$ (it contains, e.g., $3\,052\,446\,177\,238\,342\,414$). The argument above shows that the $51$ values are all of $\mathcal A_2\cap[2,0.69\cdot10^{399})$ and relates them to $\mathcal M$ and to~\eqref{eq:EM}. This method also determines A129935: since $2/(2^{1/n}-1)=2n/L-1+L/(6n)+O(n^{-3})$, the condition defining A129935 forces $0<L-2n/N<\frac{L}{3N^2}(1+O(1/n))$ with $N=\lceil2n/L\rceil$, so $2n/N$ is a convergent of $L$ (Legendre) and $N=tq_j$ with $t^2<\frac12(a_{j+1}+2)$. Testing these candidates directly (the arguments of $\lfloor\cdot\rfloor$ and $\lceil\cdot\rceil$ are all at distance more than $6\cdot10^{-400}$ from an integer, while the same kind of error analysis, given in the README accompanying \texttt{verify.py} (see Data availability), bounds the error of their computed values by $10^{-2000}$) shows that A129935 has exactly $101$ terms below $0.69\cdot10^{399}$: the $100$ terms of its b-file, the last of which has $386$ digits, and the $393$-digit element of $\mathcal A_2$ above.

\emph{A175406.} These pairs also describe this sequence, where $a(N)$ is the greatest $k$ with $(1+1/N)^k\le2$, so that $a(N)=\lfloor\ell_{N+2}\rfloor$. Since $\ell_m<L(m-\frac32)<\ell_m+1$, the formula $a(N)=\lfloor(N+\frac12)L\rfloor$ conjectured there by C.~R.~Greathouse~IV fails exactly when some integer $n$ satisfies $0<\delta(m,n)<h_m$ with $m=N+2$. It therefore fails for exactly $51$ values $N<10^{399}-2$, namely $N=m-2$ for the $51$ pairs above. The smallest, $N=1\,121\,626\,023\,352\,383$, was given as a counterexample, without a claim of minimality, by P.~W.~Alexandersson in A175406, and independently in~\cite{omega}.
\end{remark}

\begin{proposition}\label{prop:inf}
Let $j$ be even with $p_j$ even, $q_j$ odd and $q_j\ge2\cdot10^5$. Put $\eta_j=q_j(q_jL-p_j)$ and, for odd $t\ge1$, $m_t=(tq_j+3)/2$ and $n_t=tp_j/2$. Let
\[
\Theta(q)=\frac{4L\,h_{(q+3)/2}}{1+3/q}\qquad(q\ge3\text{ odd}).
\]
Then the following are equivalent:
\begin{enumerate}
\item[(a)] $(m_t,n_t)\in\mathcal M$ for some odd $t\ge1$;
\item[(b)] $0<\delta(m_t,n_t)<h_{m_t}$ for some odd $t\ge1$;
\item[(c)] $\eta_j<\Theta(q_j)$.
\end{enumerate}
Moreover $\Theta(q)=\frac L3+O(q^{-2})$.
\end{proposition}

\begin{proof}
Write $\delta_t=\delta(m_t,n_t)$. By the proof of Proposition~\ref{prop:M} (with $q=tq_j$ and $2n=tp_j$), $\delta_t=\frac{t^2\eta_j}{4L}\big(1+\frac3{tq_j}\big)=\frac{\eta_j}{4L}\big(t^2+\frac{3t}{q_j}\big)$, so that
\[
\frac{\delta_t}{\delta_1}=\frac{t^2+3t/q_j}{1+3/q_j}.
\]
Since $t\ge1$, we have $t^2\le t^2+3t/q_j\le t^2(1+3/q_j)$, and therefore
\begin{equation}\label{eq:dt}
\frac{t^2\delta_1}{1+3/q_j}\le\delta_t\le t^2\delta_1 ,
\end{equation}
and by the definition of $\Theta$, (c) says exactly that $\delta_1<h_{m_1}$. Note $\eta_j>0$ since $j$ is even.

(a)$\Rightarrow$(b) is clear, since $\delta_t>0$ because $\eta_j>0$.

(b)$\Rightarrow$(c). Suppose $0<\delta_t<h_{m_t}$ for some odd $t$. The ratio $\delta_t/\delta_1$ above is at least $1$ (its numerator minus its denominator is $(t^2-1)+3(t-1)/q_j\ge0$), so $\delta_1\le\delta_t$. Since $m_t\ge m_1$ and $h_m$ is decreasing, $h_{m_t}\le h_{m_1}$. Hence $\delta_1\le\delta_t<h_{m_t}\le h_{m_1}$, which is (c).

(c)$\Rightarrow$(a). Suppose $\delta_1<h_{m_1}$. If $\delta_1>\gamma_{m_1}$, then $t=1$ works. Otherwise $\delta_1\le\gamma_{m_1}$. Put $\Gamma=(\gz+116\cdot10^{-5})(1+3/q_j)$; since $\gz<0.0038918$ and $q_j\ge2\cdot10^5$, we have $\Gamma<0.0050518\cdot1.000015<0.00506$. Since $m_t\ge m_1\ge10^5$ for every $t$, \eqref{eq:thm55} gives $\gamma_{m_t}<\gz+116/m_t\le\gz+116\cdot10^{-5}$, hence $\gamma_{m_t}(1+3/q_j)<\Gamma$ for every odd $t$; in particular $\delta_1\le\gamma_{m_1}<\Gamma$. If $t$ satisfies $\Gamma<t^2\delta_1<\frac1{12}$, then by~\eqref{eq:dt}
\[
\gamma_{m_t}<\frac{t^2\delta_1}{1+3/q_j}\le\delta_t\le t^2\delta_1<\frac1{12}<h_{m_t},
\]
so $(m_t,n_t)\in\mathcal M$. The condition $\Gamma<t^2\delta_1<\frac1{12}$ says that $t$ lies in the open interval $(s,\rho s)$, where $s=\sqrt{\Gamma/\delta_1}$ and $\rho=(12\Gamma)^{-1/2}$. Here $s>1$ because $\delta_1<\Gamma$, and $\rho>(12\cdot0.00506)^{-1/2}=4.058\ldots>4.05$. It remains to find an odd integer in $(s,\rho s)$.
\begin{itemize}
\item If $1<s<3$, then $3\in(s,\rho s)$, since $\rho s>4.05>3$.
\item If $s\ge3$, the interval has length $(\rho-1)s>3.05\cdot3>2$, so it contains two consecutive integers, one of which is odd.
\end{itemize}
In both cases there is an odd $t\ge3$ with $(m_t,n_t)\in\mathcal M$.

Finally, the expansion of $h_m$ above, with $|G_3|=\frac1{24}$, gives $h_m=\frac1{12}+\frac1{8m}+O(m^{-2})$, and $1+3/q=1+\frac3{2m}+O(m^{-2})$ for $q=2m-3$. Hence $\Theta(q)=\frac L3+O(q^{-2})$.
\end{proof}

The proof uses the upper bound in~\eqref{eq:thm55} only in the form $\gamma_m<0.0051$ for $m\ge10^5$, and does not use the lower bound in~\eqref{eq:thm55}. (Corollary~\ref{cor:inf} below also uses $\gamma_m>0$, through Proposition~\ref{prop:M}.) Any explicit bound $\gamma_m<c$ for all $m\ge m_0$ with $c(1+\frac3{2m_0})<\frac1{108}$ would do, after replacing $2\cdot10^5$ by $2m_0$; this condition is what makes $\rho>3$. The bound $\gamma_m<1/(8L)$ of~\cite[Theorem~6.6]{MR} is not sufficient for this argument. Carrying the expansion of $h_m$ one step further ($h_m=\frac1{12}+\frac u8+\frac{49}{720}u^2+O(u^3)$ with $u=1/(m-1)$), one finds $\Theta(q)=\frac L3+\frac{4L}{45q^2}+O(q^{-3})$; here $4L/45=0.06161\ldots$, in agreement with \texttt{verify.py}.

\begin{corollary}\label{cor:inf}
The following are equivalent:
\begin{enumerate}
\item[(i)] $\mathcal M$ is infinite, that is, part~(ii) of Erd\H{o}s's conjecture holds;
\item[(ii)] O'Bryant's set $\mathcal A_2$ is infinite;
\item[(iii)] there are infinitely many even $j$ with $p_j$ even, $q_j$ odd and $\eta_j<\Theta(q_j)$.
\end{enumerate}
In particular, (i) holds if there are infinitely many even $j$ with $p_j$ even, $q_j$ odd and $a_{j+1}\ge5$. Conversely, if (i) holds, then for every $\varepsilon>0$ there are infinitely many even $j$ with $p_j$ even, $q_j$ odd and $\eta_j<\frac L3+\varepsilon$; in particular, infinitely many of these have $a_{j+1}\ge3$.
\end{corollary}

\begin{proof}
By Proposition~\ref{prop:M}, every element of $\mathcal M$ is $(m_t,n_t)$ for some even $j$ with $p_j$ even and $q_j$ odd, and some odd $t$. By Remarks~\ref{rem:oeis} and~\ref{rem:A2}, an integer $n\ge2$ lies in $\mathcal A_2$ exactly when $0<\delta(M,n)<h_M$ with $M=\lceil n/L+\frac32\rceil$. The proof of Proposition~\ref{prop:M} uses only $0<\delta<h_3$, so every such pair $(M,n)$ is also of the form $(m_t,n_t)$. Each $j$ gives only finitely many pairs, since $t^2<4Lh_3(a_{j+1}+2)$ as in the proof of Theorem~\ref{thm:M}. Hence (i), respectively (ii), holds if and only if infinitely many $j$ satisfy (a), respectively (b), of Proposition~\ref{prop:inf}. By that proposition, these are the $j$ satisfying (c), apart from the finitely many with $q_j<2\cdot10^5$. For the last two statements, recall from Lemma~\ref{lem:cf} that $\eta_j=(\alpha_{j+1}+q_{j-1}/q_j)^{-1}$ with $a_{j+1}\le\alpha_{j+1}<a_{j+1}+1$ and $0<q_{j-1}/q_j<1$, so $1/(a_{j+1}+2)<\eta_j<1/a_{j+1}$. If $a_{j+1}\ge5$, then $\eta_j<\frac15<\Theta(q_j)$ as soon as $q_j$ is large, since $\Theta(q)\to L/3=0.2310\ldots$; this gives (iii). Conversely, given (iii) and $\varepsilon>0$, all but finitely many of the $j$ in (iii) satisfy $\eta_j<\Theta(q_j)<\frac L3+\varepsilon$. If $\varepsilon<0.0002$, then $1/(a_{j+1}+2)<\eta_j<0.2313$ gives $a_{j+1}+2>4.32$, that is, $a_{j+1}\ge3$.
\end{proof}

The equivalence (ii)$\Leftrightarrow$(iii) is closely related to O'Bryant's criteria~\cite[Lemmas~11 and~12]{OB}. For $1<\theta<e^3$ he gives a necessary condition for $n\in\mathcal A_\theta$ in terms of the continued fraction of $2/\log\theta$~\cite[Lemma~11]{OB}, and he applies it to $\theta=2$ in~\cite[\S6]{OB}. For $1<\theta<e^6$ he gives a condition, in terms of the continued fractions of both $2/\log\theta$ and $1/\log\theta$, for a continuant of $2/\log\theta$ to lie in $\mathcal Q_\theta\cup\mathcal A_\theta$, where $\mathcal Q_\theta$ is an auxiliary set of good denominators for $1/\log\theta$~\cite[Lemma~12]{OB}. Corollary~\ref{cor:inf} uses only the continued fraction of $L$. The parity condition on $(p_j,q_j)$ takes the place of the set $\mathcal Q_\theta$, and the threshold $\Theta(q_j)$ is exact. A first-order form of the necessary direction, stated in terms of $1/\log\theta$, also appears in Seff's thesis~\cite{Seff}: by~\cite[Corollary~3.1.1b]{Seff}, every $n\in\mathcal A_\theta$ gives a reduced fraction $r/(2n)$ with $r$ odd and $0<r/(2n)-1/\log\theta<\log\theta/(12n^2)$. For $\theta=2$ this is condition (iii) with $\Theta(q_j)$ replaced by $L/3$. Note that (ii)$\Leftrightarrow$(iii), that is (b)$\Leftrightarrow$(c) of Proposition~\ref{prop:inf}, is elementary and uses no bound for $\gamma_m$. What is new is the equivalence (i)$\Leftrightarrow$(ii), which rests on~\eqref{eq:thm55} through Proposition~\ref{prop:inf}: part~(ii) of Erd\H{o}s's conjecture is equivalent to the infinitude of $\mathcal A_2$. Since $2$ is algebraic, this is a special case of O'Bryant's question~\cite[\S7, Problem~(5)]{OB} whether $\mathcal A_\theta$ can be proved finite, or infinite, for some algebraic $\theta$. Seff~\cite[\S\S3.3--3.4]{Seff} shows that $\mathcal A_\theta$ is infinite whenever $\log\theta>6$ is irrational, which covers every algebraic $\theta>e^6$. His examples with $\mathcal A_\theta$ finite have $\log\theta$ algebraic, so $\theta$ itself is transcendental by the Hermite--Lindemann theorem. For algebraic $\theta$ with $\log\theta<3$, and in particular for $\theta=2$, the question appears to be open.

For $q_j<10^{400}$, exactly $45$ even indices $j$ with $p_j$ even, $q_j$ odd and $q_j\ge2\cdot10^5$ satisfy $\eta_j<\Theta(q_j)$. They are exactly the indices of the $51$ elements of $\mathcal A_2$ found in the proof of Theorem~\ref{thm:M}, and each of them gives at least one element of $\mathcal M$; for $j=658$ this element has $t\ge3$. For all $132$ even indices with $p_j$ even, $q_j$ odd and $2\cdot10^5\le q_j<10^{400}$, one has $|\eta_j-\Theta(q_j)|>5.5\cdot10^{-4}$; the minimum is attained at $j=180$, which passes, and the closest failures are $j=652$ and $j=640$, with $a_{j+1}=3$ and $\eta_j=0.2435\ldots$ and $0.2443\ldots$. The only such indices with $3\le q_j<2\cdot10^5$ are $j=2,6,10$. For $j=6$ and $j=10$, (c) fails by a margin above $0.2$. For $j=2$ one has exactly $\eta_2=9L-6=\Theta(3)$, the equality case $(m,n)=(3,1)$, so (c) fails with equality. All these figures are produced by \texttt{verify.py}.

\begin{remark}[Infinitude]
By Corollary~\ref{cor:inf}, part (ii) of the conjecture is a one-sided approximation property of $\log2$ with parity conditions, in which $\gamma_m$ no longer appears. This property seems out of reach, and it is not a consequence of any statement valid for all irrational numbers. Indeed, O'Bryant~\cite[Theorems~5 and~6]{OB} gives uncountably many $\theta$ with $\mathcal A_\theta$ empty and uncountably many with $\mathcal A_\theta$ infinite, and he shows that $\mathcal A_\theta=\emptyset$ for $\theta=e^{\sqrt2}$~\cite[\S6]{OB}. In terms of the convergents $p_j/q_j$ of $\log\theta=\sqrt2$ (the argument of Proposition~\ref{prop:M} applies for all large $n$, since $\sqrt2/3<\frac12$), this example fails because of the parity condition alone. Indeed, $p_j^2-2q_j^2=\pm1$ forces every $p_j$ to be odd, while every even-index convergent satisfies $q_j(q_j\sqrt2-p_j)<\sqrt2/3$, since these constants tend to $1/(2\sqrt2)$. The following count is heuristic: none of the assumptions below is known for the continued fraction of $\log2$. Suppose that the approximation constants $q_j|q_jL-p_j|$ follow the law of Bosma, Jager and Wiedijk~\cite{BJW} (density $1/\log2$ on $(0,\frac12)$), that the parity of $j$, the parity class of $(p_j,q_j)$ (a third of the indices have $p_j$ even and $q_j$ odd~\cite{JL}) and the approximation constant behave independently, and that $\log q_j\sim\pi^2j/(12\log2)$ (L\'evy). Since $m\le x$ corresponds to $q_j\le(2x-3)/t$, summing over odd $t$ ($\sum t^{-2}=\pi^2/8$) gives
\[
\#\{(m,n)\in\mathcal M: m\le x\}\sim L(3L-2)\log x=0.05506\ldots\log x .
\]
This predicts $50.6$ elements with $m<10^{399}$; Theorem~\ref{thm:M} finds $50$, which supports the heuristic but of course proves nothing about~(ii). Similarly, O'Bryant~\cite[Theorem~4]{OB} shows that $\#(\mathcal A_\theta\cap[1,n])\sim\frac{\log\theta}{12}\log n$ for almost all $\theta>1$; for $\theta=2$ and $n=0.69\cdot10^{399}$ this gives $53.0$, against the $52$ elements of $\mathcal A_2\cap[1,0.69\cdot10^{399})$ (the $51$ above and $n=1$).
\end{remark}

\begin{remark}[Other bases]\label{rem:base}
The difficulty of part~(ii) is specific to the base $2$. Fix $\vartheta>1$, put $\lambda=\log\vartheta$, and consider the pairs with $(m-2)^n<\vartheta^{-1}(m-1)^n$ and $S_n(m)>m^n/(\vartheta-1)$; for $\vartheta=2$ these form $\mathcal M$. As above, the first condition is $n>\lambda(m-\frac32-h_m/m)$. Expanding $S_n(m)/m^n=\sum_{i\ge1}(1-i/m)^n$ to first order in $1/m$, uniformly for $|n-\lambda m|=O(1)$, one finds that the second condition is $n<\lambda\big(m-a_\vartheta+O(1/m)\big)$ with
\[
a_\vartheta=\frac{\vartheta+1}{2(\vartheta-1)},
\]
as follows from $\sum_i ix^i=x/(1-x)^2$ and $\sum_i i^2x^i=x(1+x)/(1-x)^3$ at $x=1/\vartheta$. The script \texttt{verify.py} checks this numerically for $\vartheta=1.5,2,3,5$. Now $a_\vartheta=\frac32$ only for $\vartheta=2$. For $\vartheta<2$ we have $a_\vartheta>\frac32$, so for all large $m$ the upper bound $\lambda(m-a_\vartheta+O(1/m))$ for $n$ lies below the lower bound $\lambda(m-\frac32-h_m/m)$; one therefore finds only finitely many such pairs. For $\vartheta>2$ we have $a_\vartheta<\frac32$, and the admissible $n$ form an interval of length $\lambda(\frac32-a_\vartheta)+O(1/m)$. When $\lambda$ is irrational (for instance, for every algebraic $\vartheta\ne1$, by the Hermite--Lindemann theorem), the equidistribution of $\lambda m$ modulo $1$ shows that this interval contains an integer for a set of $m$ of positive lower density. Only for $\vartheta=2$ do the first-order terms cancel, so that the question reduces to the second-order comparison of $\delta$ with $\gamma_m$ and $h_m$ treated above.
\end{remark}

\section*{Data availability}
The script \texttt{verify.py}, which reproduces every computational statement of this note (Theorems~\ref{thm:main} and~\ref{thm:M}, Remarks~\ref{rem:oeis} and~\ref{rem:A2}, the numerical statements after Corollary~\ref{cor:inf}, Remark~\ref{rem:base} and the heuristic constants), its output and a README with the a priori error analysis are available as ancillary files accompanying the arXiv version of this note.

\end{document}